\UseRawInputEncoding
 \documentclass{article}
 \usepackage{amsmath}
 
\title{Non-dimensional Newton-Puiseux Expansions}

\author{C. J. Chapman, \\ Department of Mathematics, University of Keele, \\ Staffordshire ST5 5BG, UK\\ \\
H. P. Wynn, \\ Department of Statistics, \\London School of Economics, London WC2A 2AE, UK\\ \\
M. A. Atherton, \\ Department of Mechanical Engineering, \\Brunel University London, Uxbridge UB3 3PH, UK\\ \\
R. A. Bates, \\ Rolls-Royce plc, \\ Moor Lane, Derby DE24 8BJ, UK}

\begin{document}
\maketitle

\begin{abstract}
Recent results in the theory and application of Newton-Puiseux expansions, i.e.~fractional power series solutions of equations, suggest further developments within a more abstract algebraic-geometric framework, involving in particular the theory of toric varieties and ideals.  Here, we present a number of such developments, especially in relation to the equations of van der Pol, Riccati, and Schr\"{o}dinger.  Some pure mathematical concepts we are led to are Graver, Gr\"{o}bner, lattice and circuit bases, combinatorial geometry and differential algebra, and algebraic-differential equations.  Two techniques are coordinated: classical dimensional analysis (DA) in applied mathematics and science, and a polynomial and differential-polynomial formulation of asymptotic expansions, referred to here as Newton-Puiseux (NP) expansions. The latter leads to power series with rational exponents, which depend on the choice of the dominant part of the equations, often referred to as the method of "dominant balance". The paper shows, using a new approach to DA  based on toric ideals, how dimensionally homogeneous equations may best be non-dimensionalised, and then, with examples involving differential equations, uses recent work on NP expansions to show how non-dimensional parameters affect the results. Our approach finds a natural home within computational algebraic geometry, i.e.~at the interface of abstract and algorithmic mathematics.
\end{abstract}

\vspace{3mm}

\section{Introduction}
The purpose of this paper is to demonstrate that recent new results \cite{CW} in the theory of fractional power series solutions of algebraic equations may be placed in a more general setting, and that this makes possible further developments in the theory of such equations \cite{aroca} .  In particular, it is shown that ideas from algebraic geometry, especially the theory of toric varieties and ideals \cite{ABW}, provide a basis for new computational tools which draw on, and extend, currently available tools within the flourishing research area of computational algebraic geometry \cite{cox}.  Moreover, the more general framework which we identify applies also to differential equations \cite{ayad}, for example to the van der Pol \cite{dorod, momani}, Riccati \cite{JS}, and Schr\"odinger \cite{BO} equations.

The paper runs together in a straightforward way  two techniques in applied mathematics, namely those of Dimensional Analysis (DA) and rational power perturbation expansions \cite{KCN, chow}, referred to variously as  singular perturbations, the method of dominant balance, asymptotic theory and the term we shall use, Newton-Puiseux expansions (NP).  Polynomial models are a natural vehicle for combining DA and NP and they share several features.   Both techniques  lend themselves to modern computational commutative algebra and  both start with polynomials with integer exponents which, deeper inside the construction, lead to fractional exponents \cite{mac1, mac2}.

The close connection arises, we will maintain, from the idea that the equations of many physical laws are dimensionally invariant, i.e.~the equations are expressed via a sum of terms each of which has the same dimension in mass ($M$), length ($L$) and time ($T$) (or other primitives). Each term  is the product of raw terms based on physical variables and a constant term.  DA itself concerns invariance under scale transformations, which lends itself naturally to polynomial methods, as we shall see in the next section.

The starting point in the considerable literature on NP is  the Newton Polytope, which is the convex hull of the vector of exponents in the model polynomial \cite{W, kruskal, fish}. NP derives from considering a {\em dominant facet} of the polynomial, all other terms being of smaller order in a certain sense. DA and NP come together in the following consideration. If we carry out what is called  {\em non-dimensionalisation}, that is, convert both variables and constants to non-dimensional terms,  we can see how the non-dimensional constants affect the NP expansions.  A fortuitous aspect is that both DA and NP have the same conventions for handling exponents for differential forms. 

This work is based on two papers \cite{CW} and \cite{ABW}, the former of which is a recent work on  NP and can be considered as a "proof of concept"; the latter is an earlier work on DA. The classical Buckingham method for DA is available at \cite{buck} \cite{gib} \cite{tan}.  There are many classical references on asymptotics, of which two excellent examples are \cite{BO, hinch}.  These give the rules which apply to asymptotic, as opposed to convergent, expansions. Key references on algebra are \cite{cox} and \cite{sturm}. The methods of the paper apply to problems with any number of independent parameters; recent examples in the biomedical context with many such parameters, all of importance, are \cite{ard, hoss}.

\vspace{3mm}

\subsection{Dimensional analysis using toric ideals}
We begin with a summary of the key ideas of dimensional analysis as used in this paper \cite{gib, tan}.  The basic assumption is that all main variables and, as we shall see, constants, start by having dimensions. The aim, then, is to convert to an equation in which, by contrast, variables and constants are all dimensionless. 
If $z_1, \ldots, z_k$ are the dimensions, such as as mass, length and time, respectively $M$, $L$ and $T$, then  $z_1=M$, $z_2=L$, $z_3=T$. For example, if the variable $x$, which describes force, has dimensions $MLT^{-2}$, then in the $z$-notation this would give the (Laurent) monomial term $z_1z_2z_3^{-2}$. 
Dimensionless variables are also monomials  in the  variables $x_1, \dots, x_d$  which, replaced by all their $z$-monomials, reduce to unity: all the dimensions  cancel. 

In the Buckingham method the exponents of all the variables and constants are held in a $d \times k$ matrix $A= \{a_{ij}\}$, with each row indexed  by a variable and each column by a dimension \cite{buck}. Assuming $k<d$ and that  $A$ has full rank $d$ we compute a full rank integer  kernel matrix $K$. That is,  $\mbox{rank}(K) = d-k$ and $A^TK = 0$. Each column of $K$ yields a dimensional quantity (sometimes called "group") using the following rule. For each column $q$ of $K$, let $q^+$ be the positive elements, and let $q^-$ be the negative elements. Then
\begin{equation} \label{a} 
X = \frac{x^{q^+}}{x^{q^{-}}} 
\end{equation} 
is the dimensionless quantity associated with column $q$.

The approach to dimensional analysis based on toric ideals is concisely described in \cite{ABW}, and we now summarise the main features of this approach as needed here;  more detail may be found in \cite{cox, sturm, cls}.  First, we preserve the power product notation and abuse the physics by simply writing down the bank of equations formed by representing the variables $\{x_i\}$ in terms of the dimensions $\{z_j\}$:
\begin{equation} \label{b}
x_i = \prod_{j=1}^k  z_j^{a_{ij}},\; i=1, \ldots, k.
\end{equation}
Dimensionless quantities are obtained from the  {\em elimination ideal} obtained by formal algebraic elimination of  the $\{z_j\}$ from these equations. The ideal is composed of terms of the form 
\begin{equation}
\prod x^{v^+} - \prod x^{v^-},
\end{equation}
for integer exponents $v^+$ and $v^-$,

The following is a simple example; more may be found in \cite{ABW}. Let $d=5$, $k=3$  and define 
\begin{equation} \label{c}
A = 
\left[
\begin{array}{rrr}
 1 & 1 & -2\\
0 & 1 & 0 \\
0 & 1 & 1 \\
1 & 0 & 3 \\
1 & -1 & -1\\
\end{array}
\right].
\end{equation}
The toric ideal  can be written
\begin{equation} \label{d}
\langle -x_3^4x_5 + x_2^3x_4, \;\; -x_3^3x_5^2 + x_1x_4, \;\; -x_2^3x_3 x_5 + x_1x_3^2 \rangle \;.
\end{equation}
This gives three dimensionless quantities
\begin{equation} \label{e}
X_1= \frac{x_2^3 x_4}{x_3^4 x_5},\quad X_2 = \frac{x_1x_4}{x_3^3 x_5^2},\quad X_3 = \frac{x_1x_3^2}{x_2^2x_3x_5},
\end{equation}
with exponents which are candidates for the Buckingham $K$ matrix. For example, the first generator gives  a column $(0,3,-4,1,-1)^T$.  Note that whereas a matrix $K$ would have $5-3 = 2$ columns, there are three generators of the toric ideal. A simple explanation is that the toric version has the {\em saturation} property which can be found by adjoining the condition that no variable $x_i$ is allowed to be zero \cite{cox}.

The terms derivable from the columns of the Buckingham $K$-matrix method lie in the ideal, but the non-uniqueness of  $K$  (in the sense of providing a basis for the kernel) implies that there is freedom of choice. In the same way, the toric ideal is unique, but its representation in terms of generators is in general not unique. The most familiar basis is the Gr\"obner basis which  depends on the chosen {\em monomial order}, a total order of monomials which is an important input to computer algebra packages.  In \cite{ABW} we mention the "Graver" basis, which is unique but can be quite large. The alternative "circuit bases" are in a sense minimal and are used in combinatorial geometry, particularly matroid theory. The set of columns for the $K$ matrix gives the so-called "lattice" basis. In summary, we obtain the toric ideal from the lattice basis by saturation. A canonical reference for toric ideals is \cite{sturm}, Chapter 4.

\section{The variable-constant ideal}
We shall use the notation $x^{\alpha} = x_1^{\alpha_1} x_1^{ \alpha_2} \cdots x_d^{ \alpha_d}$ for a monomial in the $n$ variables $(x_1, \ldots, x_d)$ where
$\alpha = (\alpha_1, \ldots, \alpha_d)$ is a vector of non-negative integers. A polynomial is then written
\begin{equation} \label{f}
f(x) = \sum_{\alpha \in A} c_{\alpha} x^{\alpha},
\end{equation}
where the $c_{\alpha}$ are coefficient vectors and $\mathcal A$ is a list.  The set $\mathcal A$ is the basis for the construction of the Newton Polytope covered in the next section. An example is the two-dimensional quadratic polynomial
\begin{equation} \label{g}
c_{00} + c_{10} x_1 + c_{01} x_2 + c_{20} x_1^2 + c_{11} x_1 x_2 + c_{02}x_2^2
\end{equation}
(where we have suppressed commas in the indexing), for which 
\begin{equation} \label{h}
{\mathcal A} = \{(0,0),(1,0),(0,1),(2,0), (1,1),(0,2)\}.
\end{equation}

The basic equation we will consider in this section is that of an algebraic variety defined by
\begin{equation} \label{i}
f(x)  = 0,
\end{equation}
where the variables $x_1 \ldots, x_d$ are physical variables of interest and the constants  $c_{\alpha}$ may also have physical dimensions. We introduce dimensions  $z_1, \ldots, z_k$ in the sense of dimensional analysis, of which, as mentioned, the most familiar are $z_1 = M$, $z_2 = L$, $z_3 =T$. We shall also assume that $f(x)$ is dimensionally homogeneous as defined shortly. For a variable $x_i$ or constant $c_{\alpha}$ we use the  notation $[x_i]$ and $[c_a]$ to capture the exponents in the dimensions. Thus if $x_i$ is a force variable, its dimensions are $M L T^{-2}$, with exponents
\begin{equation}
[x_i] = (1,1,-2).
\end{equation}
This bracket notation can be extended to monomials and polynomials by using the exponent rule that, for $\alpha = (\alpha_1, \ldots, \alpha_d)$,
\begin{equation}
[c_{\alpha} x^{\alpha}] = [c_{\alpha}] + \sum_{i=1}^n  \alpha_i [x_i].
\end{equation}
Dimensional homogeneity of the equation $f(x)$ is the condition that every term  $c_{\alpha}x^{\alpha}$ has the same dimension. That is, for some fixed integer vector $b$, we have
\begin{equation}
[c_{\alpha} x^{\alpha}] = \beta \;\;\; \mbox{for all} \;\;\; {\alpha} \in \mathcal A,
\end{equation}
so that 
\begin{equation} \label{common}
[c_{\alpha}] = \beta -  \sum_{i=1}^n  \alpha_i [x_i].
\end{equation}

A main idea of the paper is that of dimensional reduction, sometimes called {\it non-dimensionalisation}. That is, we want to replace the equation $f(x) = 0$ with another polynomial equation in which all variables and  constants are dimensionless.  It turns out that the development of DA using toric ideals is rather natural.

The first step is to  write down the full list of power product equations:
\begin{eqnarray}
x_i & = & 
z^{[x_i]},\;\; i = 1, \ldots, n, \\
c_{\alpha} & = & z^{[c_{\alpha}]}, \;\; \alpha \in L.
\end{eqnarray}
Here the equals sign continues to be used in a  lazy fashion, but allows us to define the  {\em variable-constant ideal} as the toric  elimination ideal from this double set of equations.  Finally,  if possible,  we select dimensionless quantities from the toric set-up which separate into constant terms and variable terms, as seen in the following example.

\vspace{3mm}

{\bf Example 1}. With $d=2$ and using variables $x,y$ and letters  $a,b,c,d$ for constants, consider
\begin{equation}
f(x,y) = a y^3 + b x y^2 + cx^2y + d x^4.
\end{equation}
We use $M,L,T,$ and in our notation  $x = ML^{-3}$, $y= ML^{-1}T^{-2}$  (the dimensions  of density and pressure respectively).  We have
\begin{equation}
[x] = (1,-3,0),\;\; [y] = (1,-1,-2).
\end{equation}
The common dimension $\beta$ is defined only up to a common multiplier.  Another way of saying this is that for a dimensionally homogeneous equation, we may, physical interpretation permitting, and ignoring zeros, divide out by the common dimensions of the individual terms. 

We now have a full set of power products, where we have used (\ref{common}) to extract the dimensions of the constants:
\begin{equation}
\begin{array}{ccl}
x & = & ML^{-3 },\\
y & = & ML^{-1}T^{-2}, \\
a & = & M^{-3}L^3 T^ 6, \\
b & = & M^{-3} L ^5 T^4, \\
c & = & M^{-3} L^7 T^2, \\
d & = & M^{-4} L^{12}.
 \end{array}
\end{equation}
The toric ideal for $x$, $y$, $a$, $b$, $c$, $d$ is the elimination ideal obtained by (algebraically) eliminating $M,L,T$ from these equations.  It is useful to count variables. There are six equations and three dimensional variables so that, after elimination,  we can hope for two dimensionless quantities in $x$, $y$, leaving one dimensionless quantity for the parameters $a$, $b$, $c$, $d$.   As mentioned, the generators of the ideal are dependent on the choice of monomial ordering for the  Gr\"obner computation. We settle on one solution which,  as required, keeps a single dimensionless quantity for the parameters and one for each original variable. They are
\begin{equation}
\begin{array}{ccl}
X & = & \frac{\displaystyle a^2d}{\displaystyle b^3} \;x, \\
&& \\
Y & = & \frac{\displaystyle a^3 d}{\displaystyle b^4} \;y, \\
&& \\
R & = & \frac{\displaystyle ac}{\displaystyle b^2} \, .
\end{array}
\end{equation}
Substituting $x$ and $y$  back into in $f(x, y)$  leads to a full dimensional reduction of the original equation, preserving the monomial stucture:
\begin{equation}
g(X,Y) = Y^3 +XY^2 + R X^2 Y + X^4 = 0.
\end{equation}
The equality signs are now no longer lazy in that we have fixed the definition of the dimensionless quantities in terms of the original variables and constants, which we assume are defined with appropriate units.

\section{Dimensional Analysis in the differential equation case}
The rules for  dimensional analysis in the differential case are not simply a convention but an important tool of science \cite{bar, guy}. Using the notation of Section 1 we have for ordinary differentials
\begin{equation}
\left[\frac{ d^s y}{d x^s}\right] = [y] - s[x],
\end{equation}
and for combined monomials and differentials:
\begin{equation}
\left[\frac{ d^s y}{d x^s} x^{\alpha}y^{\beta}\right]= (-s+\alpha) [x] + (1+\beta)[y] .
\end{equation}
From here on the toric development is the same. If we arrive at, say, two dimensionless quantities $X$ and $Y$, corresponding to $x$ and $y$ respectively, which may be partly or wholly derived from differentials, then we can substitute back for $x$ and $y$ to obtain, for example,  $d^sY/dX^s$ as a well-defined and dimensionless quantity.    

\vspace{3mm}
{\bf Example 2}: A Riccati equation \cite{JS, BO}.
The following example is a special case of the more general Riccati equation. We take 
\begin{equation}
\frac{dy}{dx} - a y^2 +  bxy^2 -  cxy =0,
\end{equation} 
with  $y$ as length (L) and $x$ as time (T). Then by forcing the equation to be dimensionally homogeneous with dimension $ML^{-1}$, as for the first term, we obtain
\begin{eqnarray}
x & = & T,\\
y & = & L ,\\
a & = & L^{-1} T^{-1}, \\
b & = & L^{-1} T^{-2},\\
c & = & T^{-2}. 
\end{eqnarray}
This gives the variable-constant ideal
\begin{equation}
\langle a^2c-b^2, \; by-c, \; a y-b, \;\; -ay+cx, \; bx-a, \; axy-1 \rangle,
\end{equation}
from which we select the invariants
\begin{eqnarray}
X & = & \frac{bx}{a}, \\
Y & = & \frac{a^2y}{b},\\
R & = & \frac{a^2c}{b^2}. 
\end{eqnarray}
Substituting for $x$ and $y$ into the original equation gives
\begin{equation}
\frac{b^2}{a^3} \frac{dY}{dX} - \frac{b^2}{a^3} Y^2 + \frac{b^2}{a^3} XY^2+   \frac{c}{a} XY = 0.
\end{equation}
Multiplying by $a^3/b^2$ we obtain the dimensionally invariant and homogeneous form with the single dimensionless parameter $R$:
\begin{equation}
\frac{dY}{dX} - (1-X)Y^2  - R XY =0.
\end{equation}
\vspace{3mm}
We shall continue this example in Section 7.

\vspace{3mm}

{\bf Example 3}. The Schr\"odinger equation \cite{BO}.
The one-dimensional Schr\"odinger equation is
\begin{equation}
-\frac{\hbar}{2 m} \frac{d^2 \psi(x)}{dx^2} + \frac{1}{2} m \omega ^2 x^2 \psi(x) = E \psi(x),
\end{equation}
with the usual definitions of variables and constants. We use the informal equals sign to give the dimensions
of each term as
\begin{eqnarray}
m & = & M,\\
 x & = &  L ,\\ 
\hbar & = & M L^2 T^{-1}, \\
 \omega & = & T^{-1}, \\
E & = & ML^2 T^{-2}. 
\end{eqnarray} 
The rule discussed for the differentials leads to the claim that every term  has dimensions of $\psi(x)$ times energy (which has dimensions: $ML^2 T^{-2}$). The conventionally accepted dimensions of $\psi(x)$ for a one-dimensional problem are $L^{-1/2}$, but we do not yet include this in the analysis.

The raw equations, set out in ideal form, are 
\begin{equation}
\langle m - M, \; x - L, \; Tw - 1, \; ET^2 - L^2M, \; -L^2m + T\hbar \rangle,
\end{equation}
giving the Gr\"obner basis 
\begin{equation}
\langle m\omega x^2 - \hbar, \; - \hbar \omega + E, \; -x + L, \; -mx^2 + T \hbar, \; T\omega - 1, \; -m + M \rangle.
\end{equation} 
The first two terms can be considered as the variable-constant elimination toric ideal  and give
two dimensionless quantities. A little experimentation was required to represent  the 
wave function $\psi$ and the energy $E$ via different toric terms. 

This gives the so-called  nonpolarisation version, as  follows. Construct the two dimensionless quantities
\begin{equation}
X = x \sqrt{\frac{m \omega}{h}},\; \tilde{E}  =  \frac{E}{\hbar \omega}. 
\end{equation}
Then, defining $\Psi(X) = \psi(x)$, we arrive at the familiar non-dimensionalised Schr\"odinger equation:
\begin{equation}\label{schr}
- \frac{d^2 \Psi(X)}{dX^2}  + X^2 \Psi(X) = 2 \tilde{E} \Psi(X).
\end{equation}

\section{Newton-Puiseux (NP) expansions in the polynomial case}
The second part of this paper gives an explanation of Newton-Puiseux expansions with a view to applying them to dimensionless equations derived via the method of the previous sections.  The development follows \cite{CW}, which presents explicit formulae for such expansions, extending the classical theory \cite{aroca, ayad, W}.  A Newton-Puiseux expansion is a fractional power-series solution of an algebraic equation or a system of such equations \cite{W, kruskal}.  A key feature of such an expansion is that its leading fractional power is determined by a geometrical object called the Newton polytope, as will be explained below.  The polytope exists in exponent space, and its vertices must be determined as an essential first step for the expansion to `get started'.

Thus our starting point is the generic polynomial variety $f(x) = \sum_{a \in \mathcal A} c_{\alpha} x^{\alpha}$, where $x$ is $d$-dimensional but we first ignore the dependence of the coefficients $c_{\alpha}$ on parameters. The  Newton polytope of $f(x)$ is the convex hull of the set $\mathcal A$: 
\begin{equation}
\mathcal C (\mathcal A) = \mbox{conv} \{\alpha: \alpha \in A\}. 
\end{equation}
In general,  an $\alpha$ in $\mathcal A$ may lie in the interior  of $\mathcal A$,  or in some lower dimensional facet, or its relative interior.

In  \cite{CW} the authors  select two $x$-variables which may be labelled  input and output variables respectively: $x_1$ and $x_d$. The object is to expand $x_d$ in special, and possible fractional,  powers of  $x_1$, which we then apply in conditions in which $x_1$ becomes small. The intermediate variables $x_2, \ldots, x_{d-1}$ play an ancillary role in determining the trajectory of  $x_1$ as it becomes small in relation to $x_d$.
 We only give here  the case in which $\mathcal A$ has a distinguished $d-1$ dimensional facet  $F$ with just a single point $\beta$ not lying in $F$ (so that this point is a vertex of $\mathcal C(\mathcal A)$).  Thus we define 
$ {\mathcal F} = {\mathcal A} \setminus \beta $ and 
\begin{equation}
F = \mbox{conv}(\mathcal F).
\end{equation}
In this case we can express  the function $f(x)$ as the sum of two parts, for which one part corresponds to the facet $F$ and the other to the isolated vertex. That is,
\begin{equation}
f(x) = \sum_{\alpha \in F} c_{\alpha} x^{\alpha} + c_{\beta} x^{\beta}
\end{equation}
for some constant $c_{\beta}$. 
We define the algebraic variety corresponding to the facet as
\begin{equation}\label{fullf}
f^*(x) =  \sum_{\alpha \in F} c_{\alpha} x^{\alpha}.
\end{equation}
The normal to the facet $F$ plays a special role which is transferred to certain identities between the $x_i \: (i=2, \ldots, d)$, and $x_1$. 

Here is  a simple example. Let $d=2$ and consider the triangle with vertices $A: (2,0),\;B: (0,2),\; C: (3,3)$. We choose the facet $\mathcal F$ to be the line between $A$ and $B$  on which $x_1+x_2 = 1$. Specifically, 
\begin{equation} 
F = \{(c_1,c_2) = c (2,0) + (1-c) (0,2) \; \mbox{for all}\;\; 0 \leq c \leq 1\}.
\end{equation} 
The corresponding set of monomials is
\begin{equation}
\{x_1^{2c} x_2^{2(1-c)}, \; 0 \leq c \leq 1\}.
\end{equation} 
Next, we ask for conditions on $x_1, x_2$ such that  the above set does not depend on $c$, or in other words what functions of $x_1,x_2$ are invariant on $F$. The condition is clearly $x_1=x_2$; that is, $x_1 x_2^{-1}$ is constant on $L$,  the value being unity.  Important for us is that the exponent vector  $(1,-1)$ is the vector orthogonal to $F$. The general case will appear in the next subsection. 

\subsection{Computation of exponents}
Following \cite{CW}, we first identify  $d$ vertices of $F$: $\alpha^{(1)}, \ldots, \alpha^{(d)}$. The vectors
\begin{equation}
V = \{\alpha^{(1)} - \alpha^{(j)}\},\; j=2, \ldots, d
\end{equation} 
form a basis of the affine subspace containing the facet $\mathcal F$. The normal $m$ to this space  (unique up to scale) is  the solution $m$ to the equations
\begin{equation} \label{orthog}
(a^{(1)} - a^{(j)} )^T m =0, \;\; j = 2, \ldots, d.
\end{equation}
We can now, under suitable general position assumptions, reconstruct a special basis for $\mbox{span}(V) $ of the form
\begin{equation}
(1, - q_2, 0, \ldots,0), \;(1,0,-q_3,0, \ldots, 0), \, \dots, \, (1, 0, \ldots, - q_d),
\end{equation}
where
\begin{equation}
q_j =   \frac{m_1}{m_j}, \;\;  j = 1,\ldots, d.
\end{equation}
Once we have decided that the basis members should have only two non-zero terms  we have uniqueness (up to scale). 

We then substitute 
$x_j = s_j x_1^{m_j/m_1} \, (j = 2, \ldots, d)$ for some choice of the $\{s_j: j = 2, \ldots, d-1\}$, leaving $s_d$ free for the moment.   
From a physical point of view, this substitution preserves dominant balance on the facet $F$ and gives a two-dimensional problem in $x_1$ and $x_d$. 

For clarity we look at a case with $d=3$. Take
\begin{equation}
f(x_1,x_2,x_3) = x_1 x_2 + x_1 x_2^2 x_3 + x_1^2 x_2 x_3 ^2 + x_1^2 x_2^2 x_3^2.
\end{equation}
The Newton polytope is
\begin{equation}
\mbox{conv} \{ (1,1,0), \; (1,2,1), \; (2,1,2), \; (2,2,2) \}.
\end{equation}
Select the convex hull of the first three points as the facet  spanned (affinely)
by 
\begin{equation}
(1,2, 1) - (1,1,0) = (0,1,1), \quad (2,1, 2) - (1,1,0) = (1,0,2). 
\end{equation}
The normal vector is  $m = (1,-2, 0)$, which leads to the subsitution $x_2 =  s_2 x_1^{1/2}$, and the vector $(0,1,1)$  (lying in the facet) yields $x_3 = s_3 x_1^{-1}$.

\subsection{The facet equation}

We have shown that, by a suitable choice of basis,  for every $ j = 2, \ldots, d$ there are a constant $s_j$ and rational exponent $m_j/m_1$ such that on the variety
$f^*(x) = 0 $ we have 
\begin{equation}
x_j = s_j x_1^{ m_j/m_1}.
\end{equation} 
In the case $j=d$ the term  $s_d x_1^{m_d/m_1}$
will appear in the first of  our NP expansions for $x_d$ in terms of $x_1$.

Substituting each $x_j$ in $f(x)$ we obtain
\begin{eqnarray}
f(x) & = &  \sum_{\alpha \in F}c_{\alpha}  \prod_{i=1}^d x_i^{\alpha_i} + c_{\beta}  \prod_{i=1}^d x_i ^{b_i} \\
      &  = & \sum_{\alpha \in F} c_a x_1^ {\alpha_1} \prod_{i=2}^d \left( s_i x_1^{ m_j/m_1}      \right)^{\alpha_i} + c_{\beta} \prod_{i=1}^d   \left( s_j x_1^{m_j/m_1}      \right) ^{b_i}.
\end{eqnarray}
Now using (\ref{orthog}) we can extract $x_1$ from the first term above to give
\begin{equation}
f(x) =  x_1^{c_1} \sum_{\alpha \in F} c_a  \prod_{i=2}^d s_d^{a_i} + c_{\beta} x_1^{c_2} \prod_{i=2}^d s_i^b,
\end{equation}
for certain constants $c_1$ and $c_2$. Again, extraction of $x_1$ is essentially by a dominant balance argument.

Now assume that $\{s_1, \ldots, s_{d-1} \}$ have been chosen  and  that  all the exponents $\{\alpha \in F\}$ and $\beta$ are integers. Then the term
\begin{equation}
F(s_d) =  \sum_{\alpha \in F} c_{\alpha} \prod_{i=2}^i s_i^{\alpha_i}
\end{equation}
is a polynomial in $s_d$. Similarly, for the last term we define
\begin{equation}
G(s_d) = -c_{\beta} x_1^{c_2} \prod_{i=2}^d s_d^{\beta}.
\end{equation}
Note that  $G(s_d)$ becomes more complex if we have more than one off-facet exponent 
(covered in \cite{CW}). 

As mentioned, we will see expansions which start with the term $s_d x_1^{m_d/m_1}$. The approach  in \cite{CW} is the following. We consider $x_d = s_dx_1^{m_d/m_1}$ to be a first approximation, where $s_d$ is a solution to the {\em facet equation} 
\begin{equation} 
F(s_d) = 0. 
\end{equation} 
Then we  seek a solution via a perturbation of the form $s_d x_1^{m_d/m_1}(1+z)$ to the full perturbation equation
\begin{equation}\label{main}
F(s_d(1+z)) = G(s_d(1+z)) x_1^c,
\end{equation}
where $c = c_2-c_1$.

\section{Implicit expansion based on $z$}
In \cite{CW} the authors give detailed formulae for computing $z$  using a version of the Sylvester method based on  the Fa\`a di Bruno formula. Here we briefly show how to set up the solution using computer algebra.

The first step  is to determine the value of $s_d$ as a chosen root of the facet equation $F(s_d) = 0$. We will assume for simplicity that this is not a multiple root.
The expansion for $z$ is a Taylor expansion in the variable $\delta = x_1^c$ and  takes the form
\begin{equation}
z = z_0 \delta ( 1 + z_1 \delta + z_2 \delta^2 + \ldots).
\end{equation}
The steps of the procedure are as follows.
\begin{enumerate}
\item Start with the initial solution $x_d^{(0)} = s_d x_1^{m_d/m_1}$ based only on a solution of the  facet equation $F(s_d)=0$.
\item Substitute the test solution $  x_d^{(0)} (1+z) $ into the main equation (\ref{main}).
\item Truncate the result of the substitution to order $O(z)$ and solve for $z$.
\item  Set $ z_1\delta $ equal to the lowest order terms in the solution of the last step.
\item  Continue the process with a new test solution $x_d^{(1)} = x_d^{(0)}(1+ z_1 \delta +z)$.
\end{enumerate}
The method  can be considered as an implicit Taylor series expansion in $\delta$, based on  $F$ and $G$. We now give an example to explain the steps (2) and (3).

{\bf Example 4} (Catalan numbers \cite{com, GKP}).  Consider the simple case with $d=2$ and  
$ x=x_1$, $y = x_2 $:
\begin{equation}
f(x,y) = y - x + x^2 y^2 = 0 .
\end{equation}
The Newton polytope is the triangle with vertices $\{ (1,0), (0,1), (2,2)\}$.  The first approximation is $y=x$. Next substitute the test solution $y= x (1+z)$ into the full equation  to obtain
\begin{equation}
-x^4 + (-2x^4 + x)z  + O(z^2). 
\end{equation}
Keep the first two terms and solve for $z$ to obtain
\begin{equation}
z_1^* =  \frac{x^3}{1-2x^3}.
\end{equation}
The lowest order terms give the solution $z_1 = x^3.$
Repeating the process, with the next trial solution  $y = x(1+x^3 +z)$, gives the series
\begin{equation}
y = x( 1+ x^3 + 2 x^6 + 5 x^9) = x + x^4 + 2x^7 + 5 x^{10} + \ldots \;. 
\end{equation} 
One of the solutions of  $f(x,y)=0$ for $y$ is
\begin{equation}
\frac{1- \sqrt{1-4 x^3}}{2 x^2},
\end{equation}
whose Taylor expansion is the same as above. Replacing the $x^3$ term by $x$ and multiplying by $x$ leads to the generating function for the Catalan numbers $1, 2, 5,\ldots$ (see \cite{com, GKP}).

\vspace{3mm}

\section{The differential equation case}
The usual approach for linear differential equations is to allow the individual terms  to contribute to an adapted Newton polytope, often referred to as the Kruskal-Newton polytope \cite{kruskal} \cite{fish}. This representation is  the same as used in DA so that the term 
\begin{equation}
\frac{d^s y}{d x^s} x^{\alpha} y^{\beta}
\end{equation}
contributes the point $(\alpha-s, \, \beta +1)$ to the Newton polytope.

The analysis then follows  the same lines as for polynomials. We select a facet $F$  of the "differential" Newton polytope and seek the analogue of the facet polynomial $F(s_d)$, having in mind that this will be a polynomial-differential equation in $x_d$. In our approach we use the following rules.
\begin{enumerate}
\item All differentials are of the form $d^s x_d/d x_1^s$.  
\item The equation is linear in the differential sense (terms such as $ \left(d^s x_d/d x_d^s \right)^2$ do not appear).
\item Assign, for $1 \leq j \leq d,$ 
\begin{equation}
x_j = s_j x_1^{m_j/m_1},
\end{equation}
with $m_j$ chosen as in the polynomial case, but using the differential polytope.
\end{enumerate}
Under these rules we have the first perturbation equation
\begin{equation}\label{maindiff}
\tilde{F}(s_d(1+z)) = \tilde{G}(s_d(1+z))x_1^c,
\end{equation}
where $\tilde{F}$ and $\tilde{G}$ are polynomial linear differential forms  in $x_1$ and $x_d$, and $s_d$ is the solution of the polynomial differential facet equation
\begin{equation}\label{diff}
\tilde{F}(s_d) = 0.
\end{equation}
The main difficulty in the differential case is that every iteration to update $z$ in the perturbation may itself require the solution of a new differential equation.  

\vspace{3mm}

{\bf Example 5}. Riccati equation \cite{JS, BO, AM}. We take the non-dimensionalised Riccati equation of Example 2 and revert to lower case, except for $R$, the dimensionless constant. This is a two-dimensional equation, but we  add a third variable $\epsilon$, taken to be small, in order to carry out a perturbation analysis \cite{KCN, hinch}. Thus the notation $(x_1,x_2,x_3)$ is more recognisable as $(\epsilon, x, y)$, and we require an NP  expansion  of $y$ in terms of $\epsilon$. The perturbed and now non-dimensionalised Riccati equation is
\begin{equation} 
\epsilon \frac{dy}{dx} - (1-x)y^2 - R x y = 0.
\end{equation}
First separate  the equation to the form
\begin{equation} 
\epsilon \frac{dy}{dx} - y^2 - R x y + x y^2 =0.
\end{equation}
Select a facet in the (differential) Newton polytope whose vertices in the space with points $(\epsilon, x, y)$ are obtained from the first three terms:
\begin{equation} 
(1,-1,1), \; (0,0,2), \; (0,1,1).
\end{equation}
The facet is spanned by the vectors
\begin{eqnarray}
(1,-1,1) - (0,1,1) & = & (1,-2,0), \\
(0,0,2) - (0,1,1) & = & (0,-1,1) .
\end{eqnarray} 
The first equation (see subsection 3.1) gives the substitution $x= s \epsilon^{1/2}$, for a fixed constant $s$. This substitution affects the differentiated term and we obtain a differential equation in $y$ and $\epsilon$ in which $s$ is taken as constant:
\begin{equation} \label{fullR}
2\frac{\epsilon^{3/2}}{s}\frac{dY}{ d \epsilon}  - Y^2 - s R \epsilon^{1/2}y   =  - s \epsilon^{1/2} y^2.
\end{equation}
The left-hand side gives a  differential equation  
\begin{equation}
2\frac{\epsilon^{3/2}}{s}\frac{dy}{ d \epsilon}  - y^2 - s R \epsilon^{1/2}y   =   0,
\end{equation}
defined by the face, and has solution
\begin{equation}
y = \frac{(1-s^2C)\sqrt{\epsilon}}{C(1-s^2 R) \epsilon^{s^2A/2}\sqrt{\epsilon} + s},
\end{equation}  
where $C$ is a constant of integration  which we will set to zero, for simplicity. This  gives our initial solution:
\begin{equation}
y_0 = \frac{(1-s^2 R) \sqrt{\epsilon}}{s}.
\end{equation}
The next trial solution is $y_1 = y_0 (1+z)$, which we substitute into the full equation (\ref{fullR}) to give a first order solution for $z$:
\begin{equation}
z_1 = - \frac{1-s^2R}{sR}.
\end{equation}
The trial solution after that is  $ z_2 = y_0 (1 + z_1 + z),$
and so on, giving
\begin{equation}
\frac{1-s^2R}{s} \sqrt{z} - \frac{(1-s^2R)^2}{s^2R} z  + \ldots.
\end{equation}
In this case the development can be verified to be the Taylor series expansion in $\sqrt{z}$ of the full solution 
\begin{equation}
\frac{R \sqrt{z}(1-s^2R )}{sR+(1-s^2R)\sqrt{z}}.
\end{equation}
Hence the series converges when
\begin{equation}
|z| < \left(\frac{sR}{1-s^2R}\right)^2.
\end{equation}
This demonstrates that convergence of the NP expansion depends on the value of the dimensionless quantity $R$ (as expected), and also on  the chosen value of the constant $s$..  

\vspace{3mm}

{\bf Example 6}. van der Pol equation \cite{dorod, momani, KCN, chow}.
There is an extensive literature on the van der Pol equation, the perturbation theory for which was closely studied as one of the founding models of dynamical systems and chaos (see in particular Section 11.4 in \cite{KCN}; Chapter 7 of \cite{chow}; and the work of Dorodnitsyn \cite{dorod}).  For an alternative approach involving modern ideas, see \cite{momani}.

The classical van der Pol differential equation is
\begin{equation}\label{original}
\frac{d^2y}{dx^2} + \mu(y^2-1)\frac{dy}{dx} + \omega^2 y = 0,
\end{equation}
where $y$ and $x$ are distance and time, respectively.

Model (\ref{original}) is already in non-dimensional form, but we wish to present it as an example of a more general version. Our starting point, then, is the equation
\begin{equation}\label{vdP}
\frac{d^2y}{dx^2} + a(y^2-y_0)\frac{dy}{dx} + b^2 y = 0,
\end{equation}
where $y_0$ is the value of $y$ chosen at a special time such as $x=0$.
We first write down the dimensions with, as before,  the constants determined so that the equation is dimensionally homogeneous. This gives
\begin{equation}
y=L,\;\; y_0 = L,\;\; x = T,\;\;  a = LT^{-1}, \;\; b^2 = T^{-2}.
\end{equation}
The ideal summarising these equations is 
\begin{equation}
\langle  y - L, \;\; y_0 - L, x - T, \;\; L^2Ta - 1, \;\; T^2b^2 - 1 \rangle,
\end{equation}
and a representation of the toric elimination ideal is
\begin{equation}
\langle  b-1, \;\; y_0^4 a^2-1, \;\; -y_0+y, \;\; -y_0^2 a+x \rangle.
\end{equation}
We can then write down the dimensional variables and constants:
\begin{equation}
X = \frac{x}{a y_0^2}, \;\; Y = \frac{Y}{y_0}, \;\; A= y_0^2 a, \;\; B=b.
\end{equation}
Substituting back into (\ref{vdP}), and cancelling a common factor, give a differential equation where all terms are dimensionless:
\begin{equation}
\frac{d^2Y}{dX^2} + A^2 (Y^2-1)\frac{dY}{dX} + A^2 B^2 Y =0.
\end{equation}

We may now assume  that the  non-dimensionalisation has already been carried out and applies to the constants $\mu$ and $\omega$. Not to complicate the notation, we revert  to lower case $x$ and $y$ . The points in the Newton polytope for $(x,y)$ are
\begin{equation}
(-2,1), \; (-1,3), \; (-1,1), \; (0,1),
\end{equation}
which lie in two, not three, dimensions.  We choose as the facet $F$ the line between the points $(0,1)$ and $(-2,1)$, which also contains the point $(-1,1)$, and then the  point $(0,1)$ is off-facet. Thus the differential equation is
\begin{eqnarray}\label{full}
\frac{d^2y}{dx^2}  -\mu \frac{dy}{dx} + \omega^2 y = -\mu y^2 \frac{dy}{dx},
\end{eqnarray}
and we  set the left-hand side to zero to obtain the facet equation. The latter has the solution
\begin{equation} 
c_1 \exp\left\{\left( \frac{\mu}{2} +\frac{\sqrt{\mu^2- 4 \omega^2}}{2} \right)x\right\}+ c_2 \exp\left\{\left( \frac{\mu}{2} -\frac{\sqrt{\mu^2- 4 \omega^2}}{2} \right)x\right\}.
\end{equation}
For simplicity we take the solution with integration constants $c_1=1$, $c_2=0$ and set
\begin{equation}
y_0 = \exp( Qx),
\end{equation}
where
\begin{equation}
Q  = \frac{\mu}{2} +\frac{\sqrt{\mu^2- 4 \omega^2}}{2}.
\end{equation}

The procedure tells us to substitute the trial solution
\begin{equation}
y = y_0 (1 + z)
\end{equation}
into the full equation (\ref{full}).
From this we obtain a differential equation for $z$. Instead of seeking a full solution by analytic methods we use inspection to test the solution $ z = C_1x$. This leads to a polynomial in $C_1$ whose constant term is
\begin{equation}
Q^2 + 2 C_1Q + \omega^2.
\end{equation}
Setting this to zero gives $C_1= - \mu/2$. The next term for $z$ is  cubic in $x$. Using the trial solution 
\begin{equation}
y = y_0\left(1- \frac{\mu}{2} x + C_2 x^3 \right)
\end{equation}
and setting the  $O(x)$ terms to  zero, we find
\begin{equation}
C_2 = \frac{\mu}{24} (-\mu^2 + 8 \omega^2 + \sqrt{ \mu - \omega^2}).
\end{equation}
Collecting gives
\begin{eqnarray}
 y \!\!\! & = & \!\!\! 1 + \frac{\sqrt{\mu^2 - 4w^2}}{2}x - \frac{w^2}{2}x^2 + \left(-\frac{\mu^3}{12} + \frac{\mu w^2}{3} - \frac{w^2 \sqrt{mu^2 - 4w^2}}{12}\right)x^3 + \ldots  \\
& = & \!\!\! 1 + \frac{V}{2} x - \frac{\omega^2}{2} x^2 - \frac{V (\mu V + V^2)}{12} x^3 + \ldots,
\end{eqnarray}
where $V =\sqrt{\mu^2- 4 \omega^2}$.

\section{Conclusion and further work}
The results in this paper are of two main types.  First, we have demonstrated that recently obtained new results giving the detailed properties of Newton-Puiseux expansions \cite{CW} may be embedded in an abstract framework involving ideas from algebraic geometry \cite{cox, sturm, cls}.  The key reason why this is a worthwhile enterprise is that algebraic geometry, although now a highly abstract subject, is at heart about the solutions of systems of polynomial equations, and therefore has much to offer of relevance to researchers in applied disciplines, where such systems of equations occur repeatedly (if not everywhere) in the modelling of complex systems \cite{cs, high, ls}.  We have demonstrated such relevance in \S2, and shown in particular, by example, how the concept of a toric ideal \cite{cls} emerges naturally in describing the way the solution of an equation scales with the coefficients.

The second type of result in the paper concerns the way in which the theory of Newton-Puiseux expansions may be combined with that of dimensional analysis to produce (i) a highly compact dominant-balance representation of the leading-order behaviour of a system; and (ii) as many correction terms as required, up to arbitrarily high order.  For both (i) and (ii) we draw extensively on the new theory developed in \cite{CW}, and by a sequence of examples involving the equations of Riccati \cite{JS, AM}, Schr\"{o}dinger \cite{BO}, and van der Pol \cite{dorod, momani}, we carry out explicitly the construction of the non-dimensional dominant balance equation.  These examples may be found in \S3 and \S6 of the paper, while \S4 gives in detail the mathematical method by which the dominant balance equation in (i) is constructed {\it ab initio}, and  \S5 indicates the key algorithmic steps (explicated in more detail in \cite{CW}) which would be required to carry out the programme (ii), if it were needed.

We believe that the above results provide a fertile basis for future developments in both the underlying theory and its use in applied science and engineering.  The modern Newton-Puiseux literature is currently in a period of explosive growth, following the fundamental papers of MacDonald \cite{mac1, mac2}, and developments in differential algebra \cite{aroca, ayad}.  Another strand is the arithmetical theory of Newton-Puiseux expansions, as expounded in \cite{bruno, kuehn}.  Two promising directions are to high-dimensional problems and partial differential equations in physics and engineering \cite{momani, ard, hoss, guy, AM, cs}, for which modern computer algebra platforms provide suitable tools; for example, special-purpose software for algebraic differential equations is increasingly included in general-purpose packages such as Maple.  Alongside this, dimensional analysis may itself be extended through the use of modern differential-algebraic methods \cite{as}.

\end{document}